\documentclass[a4paper,12pt]{amsart}
\usepackage{amssymb}
\usepackage{ifthen}
\usepackage{graphicx}
\usepackage{float}
\usepackage{caption}
\usepackage{subcaption}
\usepackage{cite}
\usepackage{amsfonts}
\usepackage{amscd}
\usepackage{amsxtra}
\usepackage{color}

\newtheorem{thm}{Theorem}
\newtheorem{cor}{Corollary}
\newtheorem{lem}{Lemma}

\newtheorem{rem}{Remark}
\newtheorem{example}{Example}
\newtheorem{defn}{Definition}
\newtheorem{prob}{Problem}

\newtheorem{conj}{Conjecture}
\theoremstyle{definition}

\newcounter {own}
\def\theown {\thesection  .\arabic{own}}

\newenvironment{pf}[1][]{%
 \vskip 3mm
 \noindent
 \ifthenelse{\equal{#1}{}}%
  {{\slshape Proof. }}%
  {{\slshape #1.} }%
 }%
{\qed\bigskip}

\newcounter{alphabet}
\newcounter{tmp}
\newenvironment{Thm}[1][]{\refstepcounter{alphabet}%
\bigskip%
\noindent%
{\bf Theorem \Alph{alphabet}}%
\ifthenelse{\equal{#1}{}}{}{ (#1)}%
{\bf .} \itshape}{\vskip 8pt}

\makeatletter
\newcommand{\Ref}[1]{\@ifundefined{r@#1}{}{\setcounter{tmp}{\ref{#1}}\Alph{tmp}}}
\makeatother

\newenvironment{Lem}[1][]{\refstepcounter{alphabet}%
\bigskip%
\noindent%
{\bf Lemma \Alph{alphabet}}%
{\bf .} \itshape}{\vskip 8pt}

\newenvironment{Conj}[1][]{\refstepcounter{alphabet}%
\bigskip%
\noindent%
{\bf Conjecture \Alph{alphabet}}%
{\bf .} \itshape}{\vskip 8pt}

\newcommand{\IC}{{\mathbb C}}
\newcommand{\ID}{{\mathbb D}}

\def\be{\begin{equation}}
\def\ee{\end{equation}}

\newcommand{\bee}{\begin{enumerate}}
\newcommand{\eee}{\end{enumerate}}

\newcommand{\blem}{\begin{lem}}
\newcommand{\elem}{\end{lem}}
\newcommand{\bthm}{\begin{thm}}
\newcommand{\ethm}{\end{thm}}
\newcommand{\bcor}{\begin{cor}}
\newcommand{\ecor}{\end{cor}}
\newcommand{\beg}{\begin{example}}
\newcommand{\eeg}{\end{example}}
\newcommand{\begs}{\begin{examples}}
\newcommand{\eegs}{\end{examples}}
\newcommand{\bdefe}{\begin{defn}}
\newcommand{\edefe}{\end{defn}}
\newcommand{\bprob}{\begin{prob}}
\newcommand{\eprob}{\end{prob}}
\newcommand{\bei}{\begin{itemize}}
\newcommand{\eei}{\end{itemize}}

\newcommand{\bcon}{\begin{conj}}
\newcommand{\econ}{\end{conj}}
\newcommand{\bcons}{\begin{conjs}}
\newcommand{\econs}{\end{conjs}}
\newcommand{\bprop}{\begin{propo}}
\newcommand{\eprop}{\end{propo}}
\newcommand{\br}{\begin{rem}}
\newcommand{\er}{\end{rem}}
\newcommand{\brs}{\begin{rems}}
\newcommand{\ers}{\end{rems}}
\newcommand{\bo}{\begin{obser}}
\newcommand{\eo}{\end{obser}}
\newcommand{\bos}{\begin{obsers}}
\newcommand{\eos}{\end{obsers}}
\newcommand{\bpf}{\begin{pf}}
\newcommand{\epf}{\end{pf}}
\newcommand{\ba}{\begin{array}}
\newcommand{\ea}{\end{array}}
\newcommand{\beq}{\begin{eqnarray}}
\newcommand{\beqq}{\begin{eqnarray*}}
\newcommand{\eeq}{\end{eqnarray}}
\newcommand{\eeqq}{\end{eqnarray*}}

\newcommand{\ds}{\displaystyle}

\newcounter{minutes}
\divide\time by 60
\newcounter{hours}
\multiply\time by 60 \addtocounter{minutes}{-\time}

\begin{document}
\bibliographystyle{amsplain}
\title[On the coefficient conjecture of Clunie and Sheil-Small on Univalent Harmonic Mappings]
{On the coefficient conjecture of Clunie and Sheil-Small on Univalent Harmonic Mappings}

\thanks{
File:~\jobname .tex,
          printed: \number\year-\number\month-\number\day,
          \thehours.\ifnum\theminutes<10{0}\fi\theminutes}

\author{S. Ponnusamy $^\dagger $
}
\address{S. Ponnusamy,
Indian Statistical Institute (ISI), Chennai Centre, SETS (Society
for Electronic Transactions and security), MGR Knowledge City, CIT
Campus, Taramani, Chennai 600 113, India. }
\email{samy@isichennai.res.in, samy@iitm.ac.in}
%
%
\author{A. Sairam Kaliraj}
\address{A. Sairam Kaliraj, Department of Mathematics,
Indian Institute of Technology Madras, Chennai--600 036, India.}
\email{sairamkaliraj@gmail.com}

\keywords{Harmonic functions, stable harmonic functions, harmonic univalent, harmonic convex,
convex in one direction, growth and covering theorems, and coefficient bound.\\
$
^\dagger$ {\tt
The first author is currently on leave from the Department of Mathematics,
Indian Institute of Technology Madras, Chennai-600 036, India.
}
}

\date{\today  
File: PonSai-5(13).tex}

\begin{abstract}
 In this paper, we first prove the coefficient conjecture of Clunie and Sheil-Small for a  class of
univalent harmonic functions which includes functions convex in some direction. Next, we prove growth and covering
theorems and some related results. Finally, we propose two conjectures. An affirmative answer to
one of which would then imply for example a solution to the conjecture of Clunie and Sheil-Small.
\end{abstract}
\thanks{ }

\maketitle
\pagestyle{myheadings}
\markboth{S.Ponnusamy, and A. Sairam Kaliraj}{Coefficient
Conjecture for Univalent Harmonic Mappings}
\section{Introduction and a Main Result}\label{PS4Sec1}
In 1984, Clunie and Sheil-Small \cite{Clunie-Small-84} proposed a conjecture on the coefficient bounds of normalized univalent
harmonic functions (see Conjecture \Ref{Har_Coeff}). This conjecture is considered to be the harmonic analog of the
Bieberbach coefficient conjecture proved by de Branges \cite{de_Branges}. The coefficient conjecture of Clunie and Sheil-Small
has been verified for a number of geometric subclasses of univalent harmonic functions but the conjecture
remains open for the full class of univalent harmonic mappings. In this article we begin the discussion by proving
the coefficient conjecture of Clunie and Sheil-Small for a larger class of univalent harmonic functions
which includes functions convex in some direction. Based on the investigation and a number of examples of this
article, we propose two new conjectures.

Let ${\mathcal H}$ denote the class of all complex-valued harmonic functions
$f=h+\overline{g}$ in the unit disk ${\mathbb D}=\{z \in {\mathbb C}:\, |z|<1\}$, where $h$ and $g$ are analytic
in $\mathbb{D}$ and normalized by $h(0)=g(0)=0=h'(0)-1 $. We call $h$ and $g$, the analytic and the co-analytic
parts of $f$, respectively, and have the following power series representation
\be\label{PSSerRep}
h(z)=z+\sum _{n=2}^{\infty}a_n z^n ~~\mbox{and}~~ g(z)=\sum _{n=1}^{\infty}b_n z^n, ~~~ z \in \mathbb{D}.
\ee
A function $f \in {\mathcal H}$ is locally univalent and sense-preserving in $\ID$ if $J_f(z)>0$ for all $z$
in $\ID$, where the Jacobian $J_f(z)$ of $f=h+\overline{g}$ is given by
$$J_f(z) = |h'(z)|^2-|g'(z)|^2.
$$
Using a result of Lewy \cite{lewy-36} and the inverse function theorem, one obtains that
$J_f(z)>0 $ in $\ID$ is a necessary and sufficient condition for $f \in {\mathcal H}$
to be locally univalent and sense-preserving in $\ID$. Consequently, $f=h+\overline{g}\in {\mathcal H}$ is
sense-preserving in $\ID$ if and only if $g'(z)=\omega (z)h'(z),$ where $\omega$ is analytic in $\ID$ with
$|\omega (z)|<1$ in $\ID$. Here $\omega$ is called the analytic dilatation of $f=h+\overline{g}$. For many basic results on
univalent harmonic mappings, we refer to the monograph of  Duren \cite{Duren:Harmonic} and also \cite{Do2, PonRasi2013}.
Denote by $\mathcal{S}_H$ the class of all sense-preserving harmonic univalent mappings $f=h+\overline{g}\in
{\mathcal H}$ and by $\mathcal{S}^0_H$ the class of functions $f \in \mathcal{S}_H$ such that
$f_{\overline{z}}(0)=0$. For the classical univalent class
$\mathcal{S} = \left\{ f=h+\overline{g} \in \mathcal{S}_H :\, g(z) \equiv 0 ~\mbox{ on }
\mathbb{D}\right\},
$
de Branges \cite{de_Branges} has proved the Bieberbach conjecture: $|a_n| \leq n$ for all $n \geq 2$.

A function $f \in \mathcal{S}_H$ is called starlike (resp. convex, close-to-convex) in $\ID$ if the range $f(\ID)$ is
starlike with respect to $0$ (resp. convex, close-to-convex), see \cite{Clunie-Small-84,Duren:Harmonic,PonRasi2013}.
The Bieberbach conjecture has been a driving force behind the development of univalent function theory so does the
coefficient conjecture of Clunie and Sheil-Small \cite{Clunie-Small-84}
for the theory of univalent harmonic mappings in the plane.

\begin{Conj}\label{Har_Coeff}{\rm \cite[Open questions]{Clunie-Small-84}}
For $f=h+\overline{g}
\in \mathcal{S}^0_H$ with the series representation as in \eqref{PSSerRep}, we have
\be\label{PSCoefConj}
\left \{
 \begin{array}{ll}
\ds |a_n| \leq \frac{(n+1)(2n+1)}{6},\\[2mm]
\ds |b_n| \leq \frac{(n-1)(2n-1)}{6}, ~\mbox{ for all }~ n \geq 2.  \\[2mm]
\ds \big| |a_n| - |b_n| \big| \leq n,
\end{array}
 \right.
\ee
\end{Conj}

The bounds are attained for the harmonic Koebe function $K(z)$, defined by
\be\label{Har_Koebe}
K(z)=\frac{z-\frac{1}{2}z^2+\frac{1}{6}z^3}{(1-z)^3}+
\overline{\left (\frac{\frac{1}{2}z^2+\frac{1}{6}z^3}{(1-z)^3}\right )}.
\ee
This conjecture has been verified for a number of subclasses of $\mathcal{S}^0_H$, namely the class of all functions
starlike, close-to-convex, convex, typically real,  or convex in one direction (see \cite{Clunie-Small-84,Sheil-Small,
WangLZhang-01}),
respectively. Recall that a domain $D\subset\mathbb{C}$ is called convex in the direction
$\theta$ $(0\leq \theta < \pi)$ if every line parallel to the line through $0$ and $e^{i\theta}$
has a connected or empty intersection with $D$. A univalent harmonic function $f$ in $\mathbb{D}$ is said
to be convex in the direction $\theta$ if $f(\mathbb{D})$ is convex in the direction $\theta$.

One of the primary issues here is to obtain useful necessary and sufficient conditions for $f$ to belong to
$\mathcal{S}^0_H$, in particular. Functions generated using such results usually have certain common properties.
For example, we have

\begin{Lem}\label{uni-theo3}{\rm \cite[Lemma 5.15]{Clunie-Small-84}}
If $h,g$ are analytic in $\mathbb D$ with $|h'(0)|>|g'(0)|$ and $h+\epsilon g$ is close-to-convex for each
$\epsilon$, $|\epsilon|=1$, then $f=h+\overline{g}$ is close-to-convex in $\ID$.
\end{Lem}

\begin{Thm}\label{lem2.1}{\rm \cite[Theorem 5.3]{Clunie-Small-84}}
A harmonic function $f=h+\overline{g} \in {\mathcal H}$ locally univalent in $\mathbb{D}$ is a
univalent mapping of $\mathbb{D}$ onto a domain convex in the
direction $\theta$ if and only if $h-e^{i2\theta}g$ is a conformal
univalent mapping of $\mathbb{D}$ onto a domain convex  in the
direction $\theta$.
\end{Thm}

For example, in order to use Theorem \Ref{lem2.1} and obtain functions that are convex in the real direction
(i.e., $\theta=0$) one adopts the following steps.
\begin{itemize}
\item Choose a conformal mapping $\phi$ with $\phi(0)=\phi'(0)-1=0$ which maps $\ID$ onto a
domain convex in the direction of the real axis.
\item Choose an analytic function $\omega: \ID \rightarrow \ID$.
\item Solve for $h$ and $g$ from $h' - g' = \phi'$ and $\omega h' - g' = 0$.
\item This gives
$$h(z)=\int_0^z \frac{\varphi'(t)}{1 - \omega(t)}\, \mathrm{d}t, ~\mbox{ and }~
g(z) = h(z) - \phi(z).
$$
\item The desired harmonic mapping is $f(z)=h(z)+\overline{g(z)}=2{\rm Re}\,(h(z))- \overline{\phi(z)}$.
\end{itemize}

For example, the harmonic Koebe function $K(z)$ defined by \eqref{Har_Koebe} is obtained by choosing $\phi(z)=z/(1-z)^2$
and $\omega(z)=z$. Similar algorithm may be formulated to construct functions that are
convex in an arbitrary direction  (see \cite{Gre}).

\begin{Thm}\label{Thm5.7}{\rm \cite[Theorem 5.7]{Clunie-Small-84}}
A harmonic function $f=h+\overline{g}\in {\mathcal H}$ locally univalent and sense-preserving in $\ID$
is convex if and only if, the analytic functions
$h(z)-e^{i\theta}g(z)$ are convex in the direction $\theta/2$ for all $\theta \in [0, 2\pi)$.
\end{Thm}

Since convex functions are convex in every direction, whenever $f=h+\overline{g}$ is convex in $\mathbb{D}$,
$h - e^{i \theta} g$ is convex in the direction $\theta / 2$ (Theorem \Ref{Thm5.7}).
More often, it is interesting to consider functions having this property. In this article we deal with
$$\mathcal{S}^0_H(\mathcal{S}) = \left\{h+\overline{g} \in \mathcal{S}^0_H :\,
h+e^{i \theta}g \in \mathcal{S}~ \mbox{for some}~~\theta \in \mathbb{R} \right\}$$
and
$$\mathcal{S}_H(\mathcal{S})=\left\{f=f_0 +
b_1 \overline{f_0}:\,f_0 \in \mathcal{S}^0_H(\mathcal{S})~~\mbox{and}~~
b_1 \in \mathbb{D} \right\}.$$ By definition, $\mathcal{S}^0_H(\mathcal{S}) \subseteq \mathcal{S}^0_H$ and
$\mathcal{S}_H(\mathcal{S}) \subseteq \mathcal{S}_H$. Moreover, it can be easily proved that
$\mathcal{S}^0_H(\mathcal{S})$ is a compact normal family. We prove that Conjecture \Ref{Har_Coeff} holds for
functions in $\mathcal{S}^0_H(\mathcal{S})$ and hence, for functions convex in one direction (see \cite{Sheil-Small}).

\bthm\label{PS5th1}
Suppose that $f=h+\overline{g} \in \mathcal{S}^0_H(\mathcal{S})$ with series representation as in \eqref{PSSerRep}.
Then \eqref{PSCoefConj} holds for all $n \geq 2$. These bounds are sharp for the class $\mathcal{S}^0_H(\mathcal{S})$.
The equality is attained for the harmonic Koebe function $K(z)$ defined by \eqref{Har_Koebe}.
\ethm

\br\label{reamark1}
If we take $h_0(z)=z+z^n$, $g_0(z)=z^n$ for $n \geq 2$ and  $F_0(z)=h_0(z)+\lambda g_0(z)$, then
$|F_0'(z)-1|<1$  in $\ID$ for any $\lambda\in \IC$ with $|\lambda +1| \leq 1/n$ and hence, $F_0(z)$ is univalent  in $\ID$ whenever
$|\lambda +1| \leq 1/n$. At the same time $f_0=h_0+\overline{g_0}$ is not locally univalent for any $n \geq 2$ (as there
are points in $\ID$ at which $h_0'(z)=0$) and hence, $f_0$ is not sense-preserving and univalent in $\ID$.
On the other hand, we do not know whether there exists at least one $\theta$ such that
$h(z)+e^{i\theta}g(z) \in \mathcal{S}$ whenever $h+\overline{g} \in \mathcal{S}^0_H$
(see Conjecture \ref{PS5conj2} at the end of Section \ref{PS5sec4}).
\er

The paper is organized as follows. We present the proof of Theorem \ref{PS5th1} in Section \ref{sec2a}.
In Section \ref{PS5sec2}, we discuss several interesting examples
of univalent harmonic functions which belongs to the class $\mathcal{S}^0_H(\mathcal{S})$. In Section \ref{PS5sec3},
we recall some important results on affine and linear invariant families of univalent harmonic functions.
As an application of Theorem \ref{PS5th1}, in Section \ref{PS5sec4}, we derive growth and covering theorems
and sharp bounds on the Jacobian and curvature of $f$ for functions
$f \in \mathcal{S}_H(\mathcal{S})$.


The present investigation together with standard examples of univalent harmonic mappings
and Theorems \ref{PS5th1} to \ref{thm_curv_sub}
(see also Remark \ref{reamark1}) suggest  the following

\begin{conj}\label{PS5conj2}
$\mathcal{S}^0_H = \mathcal{S}^0_H(\mathcal{S})$. That is, for every function $f=h+\overline{g}
\in \mathcal{S}^0_H$, there exists at least one $\theta \in \mathbb{R}$ such that $h + e^{i\theta} g \in \mathcal{S}$.
\end{conj}

It is natural to introduce and state analogous results for
$$\mathcal{C}^0_H(\mathcal{C}) = \left\{h+\overline{g} \in \mathcal{C}^0_H :\,
h+e^{i \theta}g \in \mathcal{C}~ \mbox{for some}~~\theta \in \mathbb{R} \right\}
$$
and
$$\mathcal{C}_H(\mathcal{C})=\left\{f=f_0 +
b_1 \overline{f_0}:\, f_0 \in \mathcal{C}^0_H(\mathcal{C})~~\mbox{and}~~ b_1 \in \mathbb{D} \right\}.
$$
Here $\mathcal C$ and $\mathcal{C}^0_H$ denote the class of functions $f$ from $\mathcal S$
and $\mathcal{S}^0_H$, respectively such that $f(\ID)$ is convex.
Note that $|a_n|\leq 1$ for $f\in \mathcal C$.


\section{Proof of Theorem \ref{PS5th1}}\label{sec2a}
Let $f=h+\overline{g} \in \mathcal{S}^0_H(\mathcal{S})$, where $h$ and $g$ have the power series given by
\eqref{PSSerRep}. Then $\varphi(z) = h(z) + \epsilon g(z) = z +
\sum _{n=2}^{\infty}\varphi_n z^n \in \mathcal{S}$ for some $\epsilon$ such that $|\epsilon|=1$. By the de Branges
theorem \cite{de_Branges}, $|\varphi_n| \leq n$ for all $n \geq 2$. Since $f$ is sense-preserving in
$ \mathbb{D}$, there exists an analytic function $\omega(z)$ in $\mathbb{D}$ such that $\omega(0)=0$ and
$ |\omega(z)| = |g'(z)/h'(z)|< 1 $ for all $z \in \mathbb{D}$ from which we easily obtain that $\varphi'(z)=h'(z)
(1+\epsilon \omega(z))$.
Then
\be\label{PS5eq1}
h(z) = \int_0^z \frac{\varphi'(\zeta)}{1+\epsilon \omega(\zeta)} \, \mathrm{d}\zeta ~~\mbox{ and }~~
g(z) = \int_0^z \frac{\varphi'(\zeta)\omega(\zeta)}{1+\epsilon \omega(\zeta)} \, \mathrm{d}\zeta. \nonumber
\ee
Let
$$
\frac{\omega(z)}{1+\epsilon \omega(z)} = \sum _{n=1}^{\infty}\omega_n z^n.
$$
Since $|\omega(z)| < 1$, in terms of subordination, we can write

$$
\frac{-\epsilon\omega(z)}{1+\epsilon \omega(z)} \prec \frac{z}{1-z}, ~~ z \in \ID.
$$
Here $\prec$ denotes the usual subordination (see \cite{Duren-book1}).
Since $z/(1-z)$ is convex in $\mathbb{D}$, according to the result of Rogosinski \cite{Rogosinski_PLMS} (see also
\cite[p.195, Theorem 6.4]{Duren-book1}) it follows that $|\omega_n| \leq 1$ for all $n \geq 1$.
Thus, we have
\beqq
g'(z) &=&  \frac{\varphi'(z)\omega(z)}{1+\epsilon \omega(z)} =  \left ( \varphi_1 + \sum _{n=1}^{\infty}
 (n+1)\varphi_{n+1} {z^n} \right )
 \left (\sum _{n=1}^{\infty} \omega_n z^n \right )\\
 &=&  \sum _{n=2}^{\infty} \left ( \sum _{k=0}^{n-2} (k+1)\varphi_{k+1} \omega_{n-1-k} \right) {z^{n-1}}.
\eeqq
Therefore
\beq\label{PS5eq3}
n|b_n| & \leq &  \sum _{k=0}^{n-2} (k+1)|\varphi_{k+1}| ~\mbox{ (since }~ |\omega_n| \leq 1  ~\mbox{ for all
 }~ n \geq 1) \\
  & \leq & \sum _{k=1}^{n-1} k^2 ~\mbox{ (since }~ |\varphi_n| \leq n  ~\mbox{ for all
 }~ n \geq 1) \nonumber \\
  & = & \frac{(n-1)n(2n-1)}{6}. \nonumber
\eeq
From the definition of $\varphi(z)$, we have $h(z)=\varphi(z)-\epsilon g(z)$. Therefore, one has
\be\label{PS5eq4}
|a_n| = |\varphi_n - \epsilon b_n| \leq  |\varphi_n| + |b_n|  \leq  \frac{(n+1)(2n+1)}{6}~\mbox{ for }~ n \geq 2.
\ee
This completes the proof of Theorem \ref{PS5th1}. \hfill $\Box$

\vspace{6pt}

\bcor\label{cor2}
Suppose that $f=h+\overline{g} \in \mathcal{C}_H^0(\mathcal{C})$ with the series representation as in \eqref{PSSerRep}.
Then
$$|b_n| \leq \frac{n-1}{2} ~\mbox{ and }~ |a_n| \leq \frac{n+1}{2} ~\mbox{ for all }~ n \geq 2.
$$
The bounds are attained for the half-plane mapping $f_3=h_3+\overline{g_3}$, where
$$ h_3(z)=\frac{z-z^2/2}{(1-z)^2} ~\mbox{ and }~ g_3(z)=\frac{-z^2/2}{(1-z)^2}.
$$
\ecor
\bpf
Apply the proof of Theorem \ref{PS5th1} with $|\varphi_n| \leq 1$ for all $n \geq 2$.
Then from \eqref{PS5eq3} and \eqref{PS5eq4}, the desired inequalities follow.
We remark that the function $f_3(z)$ may be obtained by choosing $\phi(z)=z/(1-z)^2$ and
$\omega(z)=-z$ in the algorithm mentioned after Theorem \Ref{lem2.1}.
The function $f_3$ is univalent and convex in $\mathbb{D}$ with $f_3(\ID)=\left\{w:\, {\rm Re}\,w>-1/2 \right\}$.
\epf

The following corollaries are easy to obtain

\bcor\label{PS5Cor1}
Suppose that $f=h+\overline{g} \in \mathcal{S}_H(\mathcal{S})$ with the series representation as in \eqref{PSSerRep}.
Then
$$|a_n| < \frac{1}{3}(2n^2+1) ~\mbox{ and }~ |b_n| < \frac{1}{3}(2n^2+1) ~\mbox{ for all }~ n \geq 2.
$$
\ecor

\bcor\label{cor3-new}
Suppose that $f=h+\overline{g} \in \mathcal{C}_H(\mathcal{C})$ with the series representation as in \eqref{PSSerRep}.
Then
$$|a_n| < n~\mbox{ and }~ |b_n| < n ~\mbox{ for all }~ n \geq 2.
$$
\ecor
We end the section with the following

\begin{conj}
Let $f=h+\overline{g} \in \mathcal{C}^0_H$, where $\mathcal{C}^0_H$ denotes the class of functions in
$\mathcal{S}^0_H$ such that $f(\ID)$ is convex. Then there exists a $\theta$ such that
the analytic function $h(z)+e^{i\theta}g(z)$ is univalent and maps $\ID$ onto a convex domain.
That is, $\mathcal{C}^0_H = \mathcal{C}^0_H(\mathcal{C})$.
\end{conj}

\section{Interesting members of the family $\mathcal{S}^0_H(\mathcal{S})$}\label{PS5sec2}
\subsection{Various Examples}
\begin{figure}[htp]
 \centering
   \begin{subfigure}[b]{0.5\textwidth}
     \centering
     \includegraphics[height=6cm, width=5.5cm]{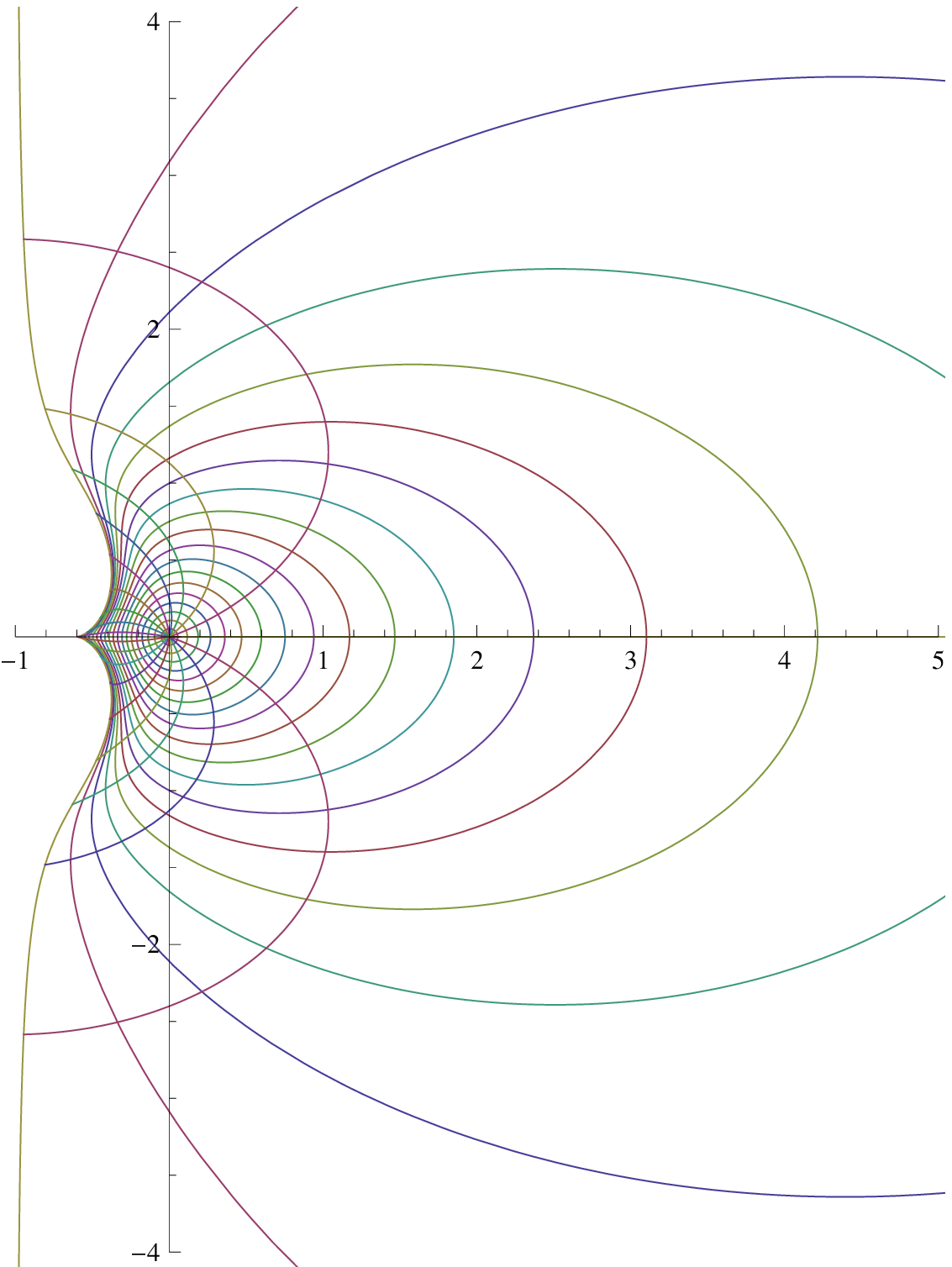}
     \caption{The range  $f_{2}(\ID)$} for $n=3$ and $\alpha=1/5$
     \label{PS4_Fig1a}
   \end{subfigure}
   ~~~
    \begin{subfigure}[b]{0.5\textwidth}
      \centering
      \includegraphics[height=6cm, width=5.5cm]{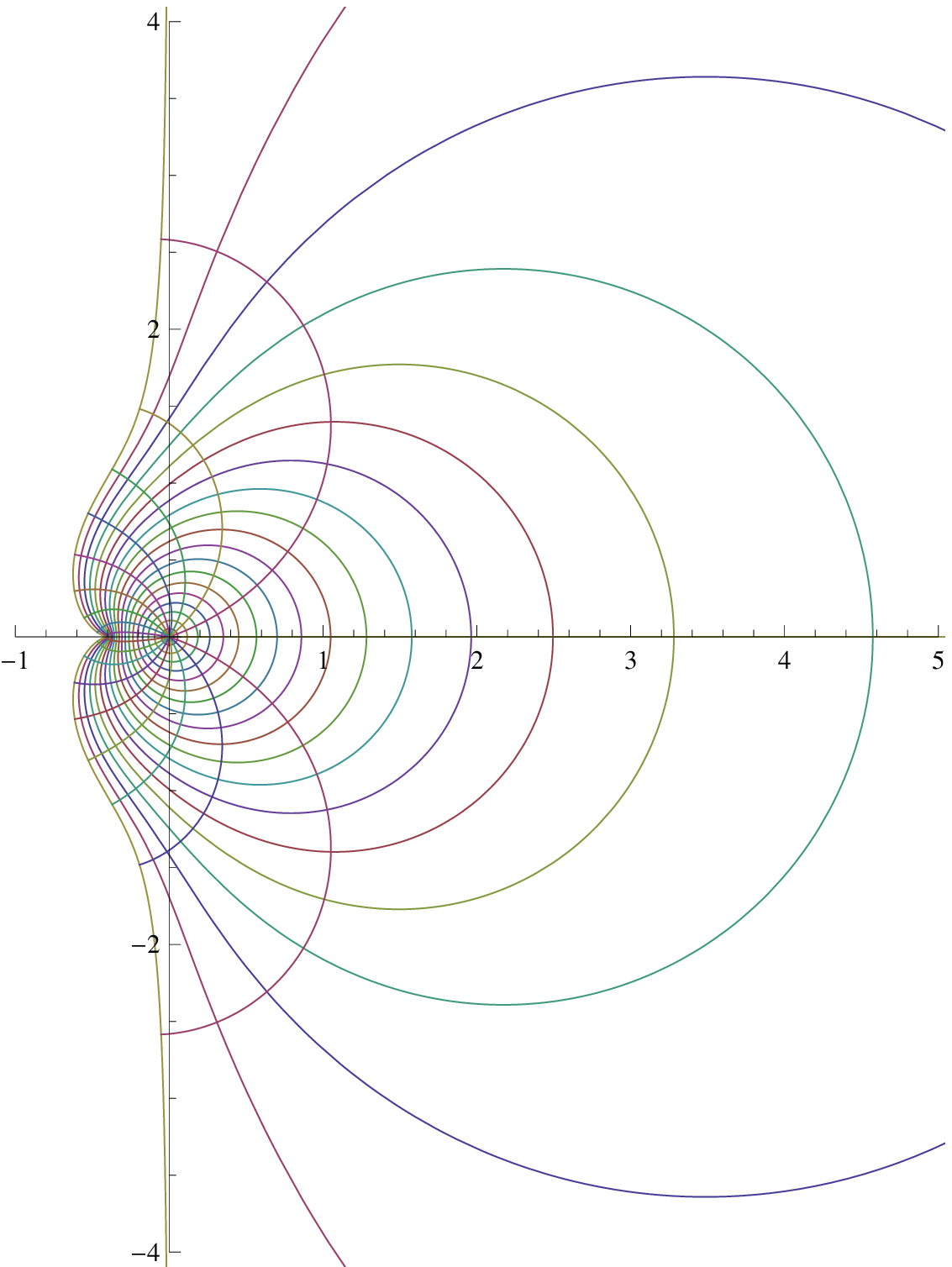}
      \caption{The range $\phi_{2, \pi}(\ID)$} for $n=3$ and $\alpha=1/5$
      \label{PS4_Fig1b}
   \end{subfigure}\\
  \begin{subfigure}[b]{0.5\textwidth}
     \centering
     \vspace{0.1cm}
     \includegraphics[height=6cm, width=5.5cm]{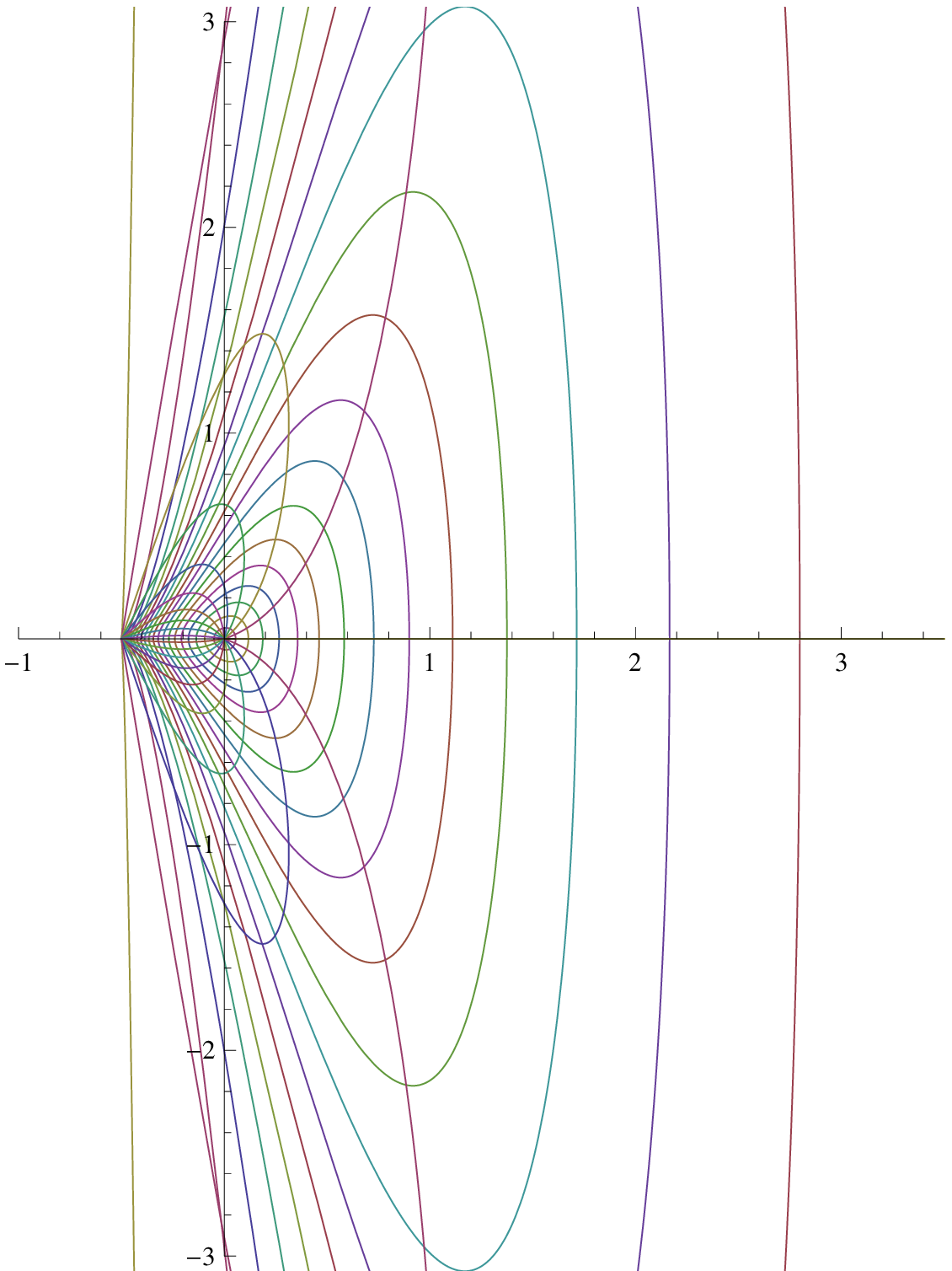}
     \caption{The range  $f_3(\ID)$}
     \label{PS4_Fig1c}
   \end{subfigure}
   ~~~
    \begin{subfigure}[b]{0.5\textwidth}
      \centering
      \includegraphics[height=6cm, width=5.5cm]{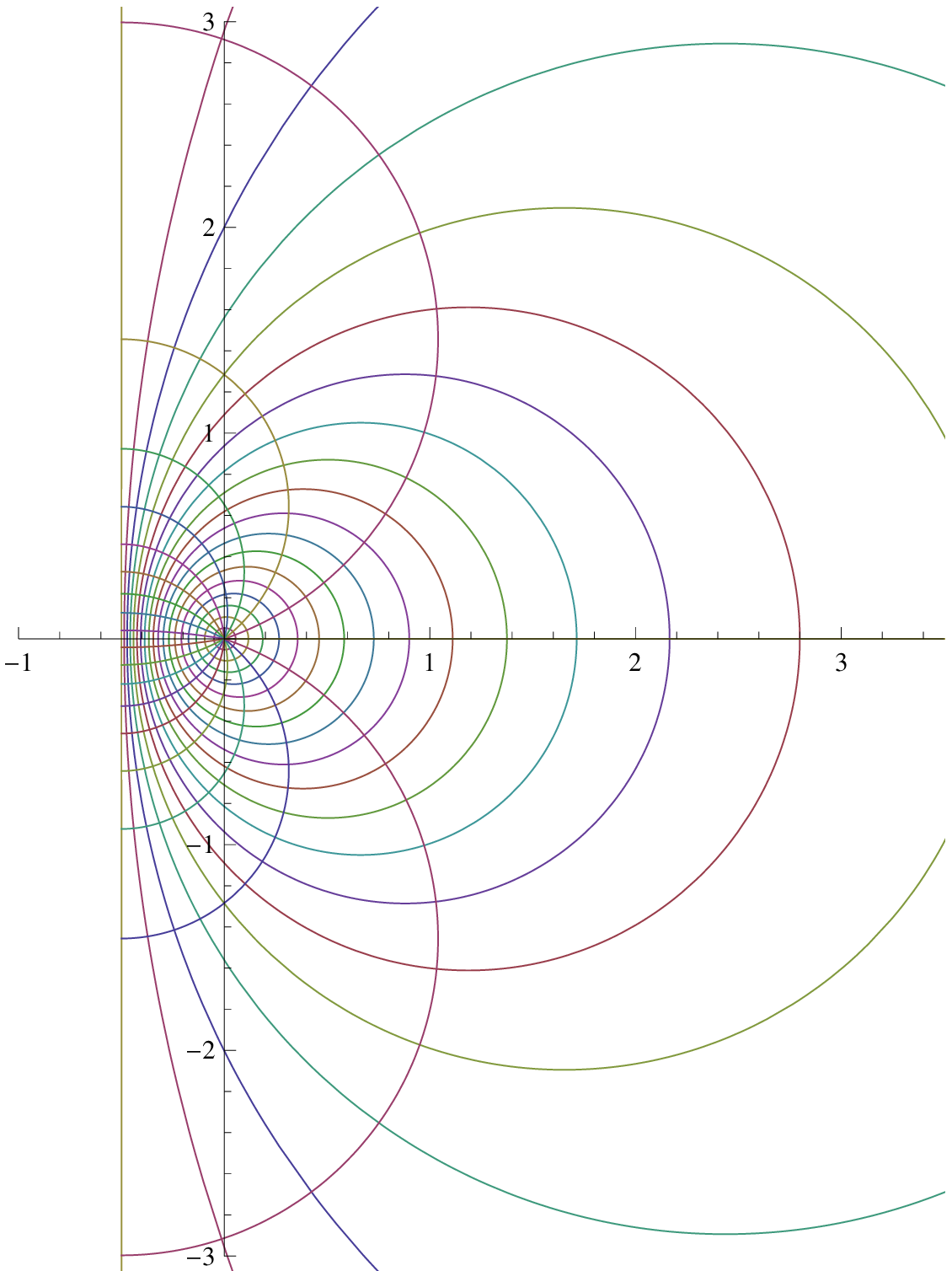}
      \caption{The range  $\phi_{3,0}(\ID)$}
      \label{PS4_Fig1d}
   \end{subfigure}\\
   \begin{subfigure}[b]{0.5\textwidth}
     \centering
     \vspace{0.1cm}
     \includegraphics[height=6cm, width=5.5cm]{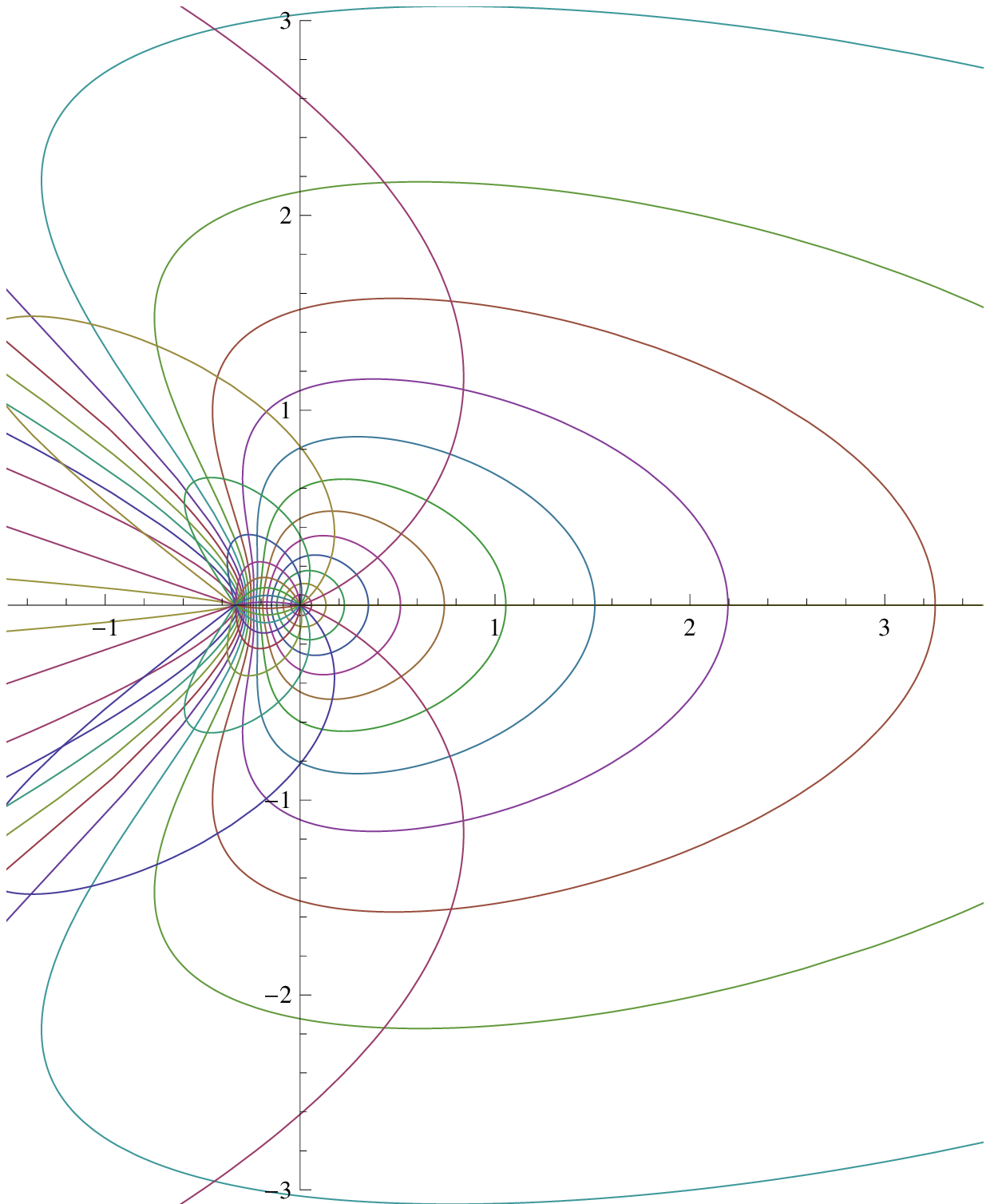}
     \caption{The range $f_4(\ID)$}
     \label{PS4_Fig1e}
   \end{subfigure}
   ~~~
   \begin{subfigure}[b]{0.5\textwidth}
      \centering
      \includegraphics[height=6cm, width=5.5cm]{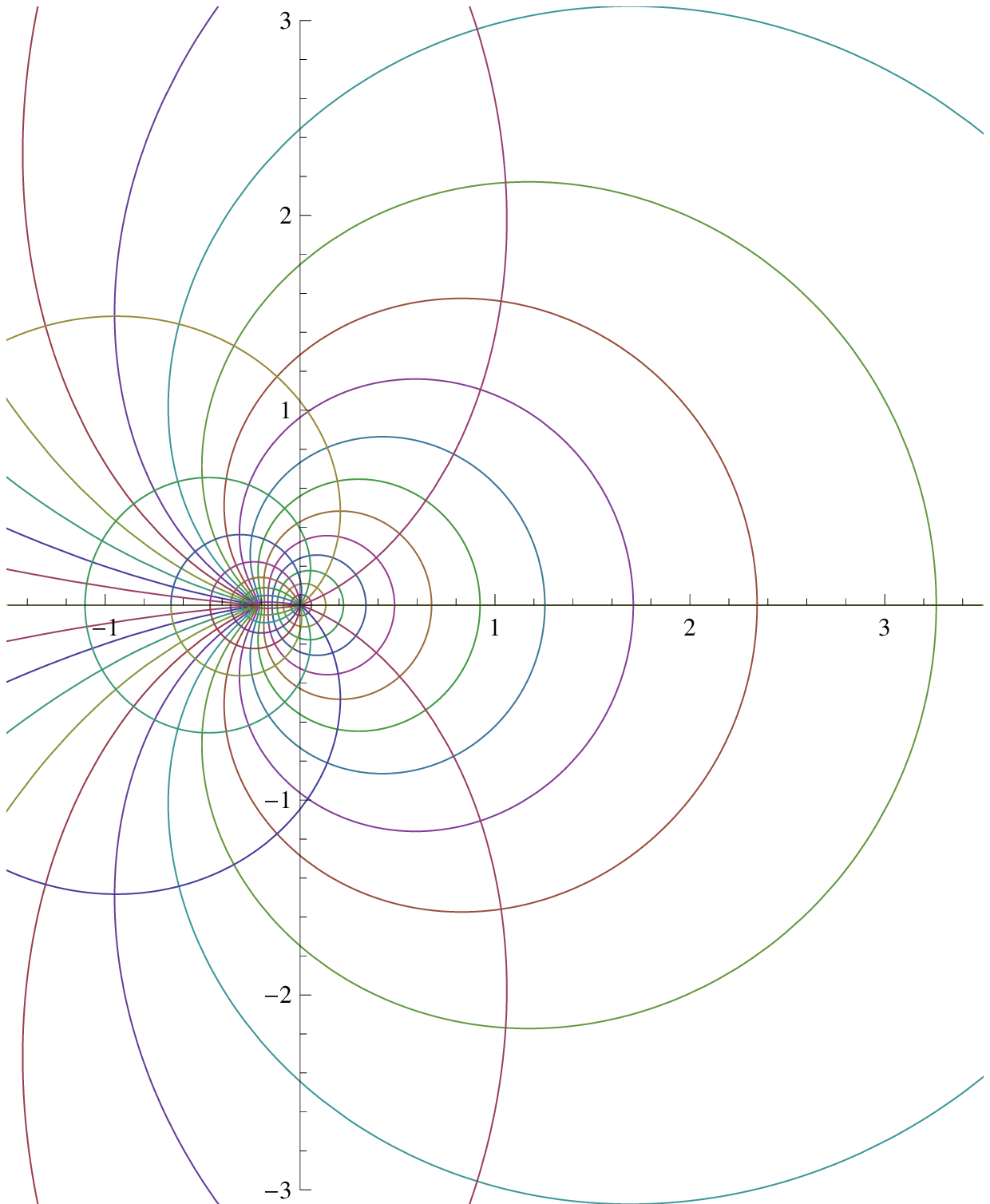}
      \caption{The range  $\phi_{4,\pi}(\ID)$}
      \label{PS4_Fig1f}
   \end{subfigure}

\caption{The images of unit disk under $f_j(z)$ and $\phi_{j,~\theta}(z)$ for $j=2,3,4$ for certain values of $\theta$.}
\label{PS5_Fig1}
\end{figure}

\begin{itemize}
\item[(1)] Let $g(z)$ be analytic in $\mathbb{D}$ with $|g'(z)| < n$ for all $z \in \mathbb{D}$ and
for some $n \in {\mathbb N}$, eg. $g(z)=z^n$, $n \geq 2$. Then the harmonic function $f_1(z) = z + \overline{g(z)/n}$
is univalent and close-to-convex in   $\mathbb{D}$. The analytic functions $\phi_{1,\theta}(z)=z+e^{i\theta}g(z)/n$ are univalent and close-to-convex for every $\theta \in \mathbb{R}$.

\item[(2)] The function $f_2(z) = z/(1-z) + \overline{g(z)}$ is univalent and close-to-convex in $\ID$,
whenever $g(z)$ is analytic in $\ID$ such that $|g'(z)|<1/|1-z|^2$ . In particular, if $\alpha \in
\IC \smallsetminus \{0\}$ such that $|\alpha| \leq 1/(2n-1)$ and $n \in {\mathbb N}$,
then the harmonic function
$$f_{2}(z) = \frac{z}{1-z} + \overline{\left (\frac{\alpha z^n}{1-z}\right )}
$$
is univalent and close-to-convex in $\mathbb{D}$. The analytic functions
$$\phi_{2,\theta}(z)=\frac{z}{1-z}+e^{i\theta}\frac{\alpha z^n}{1-z}
$$
are univalent and close-to-convex for every $\theta \in \mathbb{R}$.

\item[(3)] Consider the half-plane mapping given by $f_3$ given in Corollary \ref{cor2}, i.e.,
$$f_3(z) = \frac{1}{2}\left[ \frac{z}{1-z}+\frac{z}{(1-z)^2}\right] +\frac{1}{2} \overline{\left[
\frac{z}{1-z}-\frac{z}{(1-z)^2}\right]}.
$$
The analytic functions (see Theorem \Ref{lem2.1})
$$\phi_{3,\theta}(z)=\frac{1}{2}\left[ \frac{z}{1-z}+\frac{z}{(1-z)^2}\right] +\frac{e^{i\theta}}{2}\left[
\frac{z}{1-z}-\frac{z}{(1-z)^2}\right]
$$
are univalent and convex in the direction $\theta / 2$ for every $\theta \in \mathbb{R}$.

\item[(4)]The most interesting example is the harmonic Koebe function $K(z)=h(z)+\overline{g(z)}$
given by \eqref{Har_Koebe}. For $\theta \in [0, 2\pi)$, define
\beqq
\varphi_\theta(z)&=&\frac{z-\frac{1}{2}z^2+\frac{1}{6}z^3}{(1-z)^3}+e^{i\theta}~
\frac{\frac{1}{2}z^2+\frac{1}{6}z^3}{(1-z)^3} \\
&=&h(z)+e^{i\theta}g(z)\\
&=&z+\sum _{n=2}^{\infty}\varphi_{\theta, n}z^n,
\eeqq
where
$$\varphi_{\theta, n} = \frac{1}{6} \left( 2n^2 (1+e^{i\theta}) + 3n(1-e^{i\theta})+
(1+e^{i\theta})\right) ~\mbox{ for all }~ n \geq 2.
$$
Then $\varphi_\theta(z)$ is univalent only for $\theta = \pi$. For $\varphi_\theta$ to be
univalent in $\ID$, it is necessary that $|\varphi_{\theta,n}| \leq n$ for all $n \geq 2$. When $\theta = \pi$,
$\varphi_\theta(z)$ reduces to Koebe function $k(z)=z/(1-z)^2$, which is univalent in $\mathbb{D}$. For $\theta \in
[0, 2\pi)\smallsetminus \{\pi\}$, $|\varphi_{\theta, n}| > n$ for large values of $n$ and hence $\varphi_\theta(z)$
is not univalent in $\mathbb{D}$. Thus, $K(z)=h(z)+\overline{g(z)}$ is a member of the family $\mathcal{S}^0_H(\mathcal{S})$.

\item[(5)]As another interesting example, we consider
\beqq
f_4(z)&=&\frac{1 - (1 - z)^3}{3(1 - z)^3}+ \overline{\left (\frac{z^3}{3(1 - z)^3}\right )}\\
&=&\sum _{n=1}^{\infty}\frac{(n+1)(n+2)}{6}z^n
+\overline{\sum _{n=3}^{\infty}\frac{(n-1)(n-2)}{6}z^n}.\nonumber
\eeqq
The harmonic function $f_4(z)$ is univalent and convex in the real direction
in $\mathbb{D}$. This is because the corresponding $h-g$ is the Koebe function $k(z)=z/(1-z)^2$
and the dilatation of $f_4=h+\overline{g}$ is $\omega(z)=z^2$ (use Theorem \Ref{lem2.1} with $\theta = 0$). On the other hand,
we see that the analytic function
$$\phi_{4,\theta}(z)=\frac{1 - (1 - z)^3}{3(1 - z)^3}+ e^{i\theta}\frac{z^3}{3(1 - z)^3}
$$
is univalent only for $\theta = \pi$, and not for other values of $\theta$. The proof follows from the same reasoning as
in the previous example.
An interesting fact is that the coefficients of $f_4(z)$ are smaller than the coefficients of harmonic Koebe
function, but it does satisfy the condition $\big| |a_n| - |b_n| \big| = n$ for all $n \geq 2$
(compare with the conjecture on coefficient bounds in \eqref{PSCoefConj}).

\end{itemize}
The images of the unit disk $\ID$ under these functions for certain values of $\theta$ are shown in Figure \ref{PS5_Fig1} as plots of the images of equally spaced radial segments and concentric circles. These figures are drawn by using Mathematica.

\subsection{Stable Harmonic Univalent Functions}
Recently,  Hern\'{a}ndez and Mart\'{i}n \cite{Rodri-Maria} studied stable
harmonic univalent functions. A sense-preserving harmonic mapping $f = h +\overline{g}$ is said to be {\it stable
harmonic univalent} or simply {\bf SHU} in the unit disk (resp. {\it stable harmonic convex} ({\bf SHC}), {\it stable
harmonic starlike with respect to the origin} (${\bf SHS}^*$), or {\it stable harmonic close-to-convex}
({\bf SHCC})) if all the mappings $f_\lambda = h + \lambda g$ with $|\lambda| = 1$ are univalent
(resp. {\it convex, starlike with respect to the origin, or close-to-convex}) in $\ID$. They
proved that for all $|\lambda|=1$, the functions
$f_{\lambda}=h+\lambda\overline{g}$ are univalent (resp. {\it close-to-convex, starlike, or convex}) if and only if the
analytic functions $F_{\lambda}= h + \lambda g$ are univalent (resp. {\it close-to-convex, starlike, or convex}) for all such
$\lambda$. Let us consider
\be\label{faz}
f_{a, \lambda}(z)= a \log\left(\frac{a}{a-z}\right) +\lambda \overline{\left(a \log\left(\frac{a}{a-z}\right)
- z\right)},
\ee
where $|a|\geq 1$ and $|\lambda|=1$. A simple calculation shows that $f_{a, \lambda}$ is sense-preserving in $\ID$ and
$a \log\left(\frac{a}{a-z}\right)$ is a convex function in $\ID$. From \cite[Theorem 5.17]{Clunie-Small-84}, it
can be easily verified that $f_{a, \lambda}$ is univalent and close-to-convex in $\ID$ for all $\lambda$ such that $|\lambda|=1$.
From the above mentioned result of Hern\'{a}ndez and Mart\'{i}n \cite{Rodri-Maria}, it follows that the function
$$\phi_{a, \lambda}(z)= a \log\left(\frac{a}{a-z}\right) +\lambda \left(a \log\left(\frac{a}{a-z}\right)
- z\right),
$$
is univalent and close-to-convex in $\ID$ for all $\lambda$ such that $|\lambda|=1$. In fact in this case we
can obtain a stronger conclusion when $|a| \geq 1+\sqrt{2}$, $a \in \mathbb{R}$. With $a>1$, we compute
$$\phi'_{a, \lambda}(z)=\frac{a+\lambda z}{a-z} ~\mbox{ and }~
1 + \frac{z\phi_{a, \lambda}''(z)}{\phi_{a, \lambda}'(z)} =1+\frac{\lambda z}{a+\lambda z}+\frac{z}{a-z}.
$$
Considering the images of $|z|=r$ under $w=\frac{\lambda z}{a+\lambda z}$ and $w_1=\frac{z}{a-z}$, it follows
easily that
\beqq
{\rm Re} \left (1 + \frac{z\phi_{a, \lambda}''(z)}{\phi_{a, \lambda}'(z)} \right ) &\geq&
1 - \frac{r}{a - r}-\frac{r}{a + r}
> \frac{a^2 - 2 a  - 1}{a^2 - 1}
\eeqq
which is non-negative provided $a \geq 1+\sqrt{2}.$

This observation implies that $\phi_{a, \lambda}(z)$ is convex in $\ID$ for each $\lambda$ such that $|\lambda|=1$ and
$a\geq1+\sqrt{2}$. We see that $-\phi_{-a, \lambda}(-z)=\phi_{a, \lambda}(z)$ and hence, we conclude that
$f_{a, \lambda}$ is also a convex function for all $|\lambda|=1$, and for each
$a\in (-\infty,-1-\sqrt{2}]\cup[1+\sqrt{2}, \infty)$.

\subsection{Analog of Alexander Transform for Stable Harmonic Functions}
It is well known that a fully starlike harmonic function need not be univalent in $\ID$ (see \cite{Duren_Osgood}).
On the other hand, it is proved in \cite{Rodri-Maria} that the stable harmonic starlike functions
are necessarily univalent in $\ID$. Moreover, analog of the classical Alexander's theorem for
analytic functions has been  proved for stable harmonic functions (see \cite{Rodri-Maria}) and it reads as follows.

\begin{Thm}\label{SH_Alexander}
Assuming that the analytic functions $h$, $g$, $H$, and $G$ defined in $\ID$ are related
by
$$zh'(z) = H(z), ~\mbox{ and }~ zg'(z)=−G(z),
$$
we have that $F=H+\overline{G}$ is ${\bf SHS}^*$ if and only if $f=h+\overline{g}$ is {\bf SHC}.
\end{Thm}

\begin{figure}
 \centering
   \begin{subfigure}[b]{0.5\textwidth}
     \centering
     \includegraphics[height=6cm, width=5.5cm]{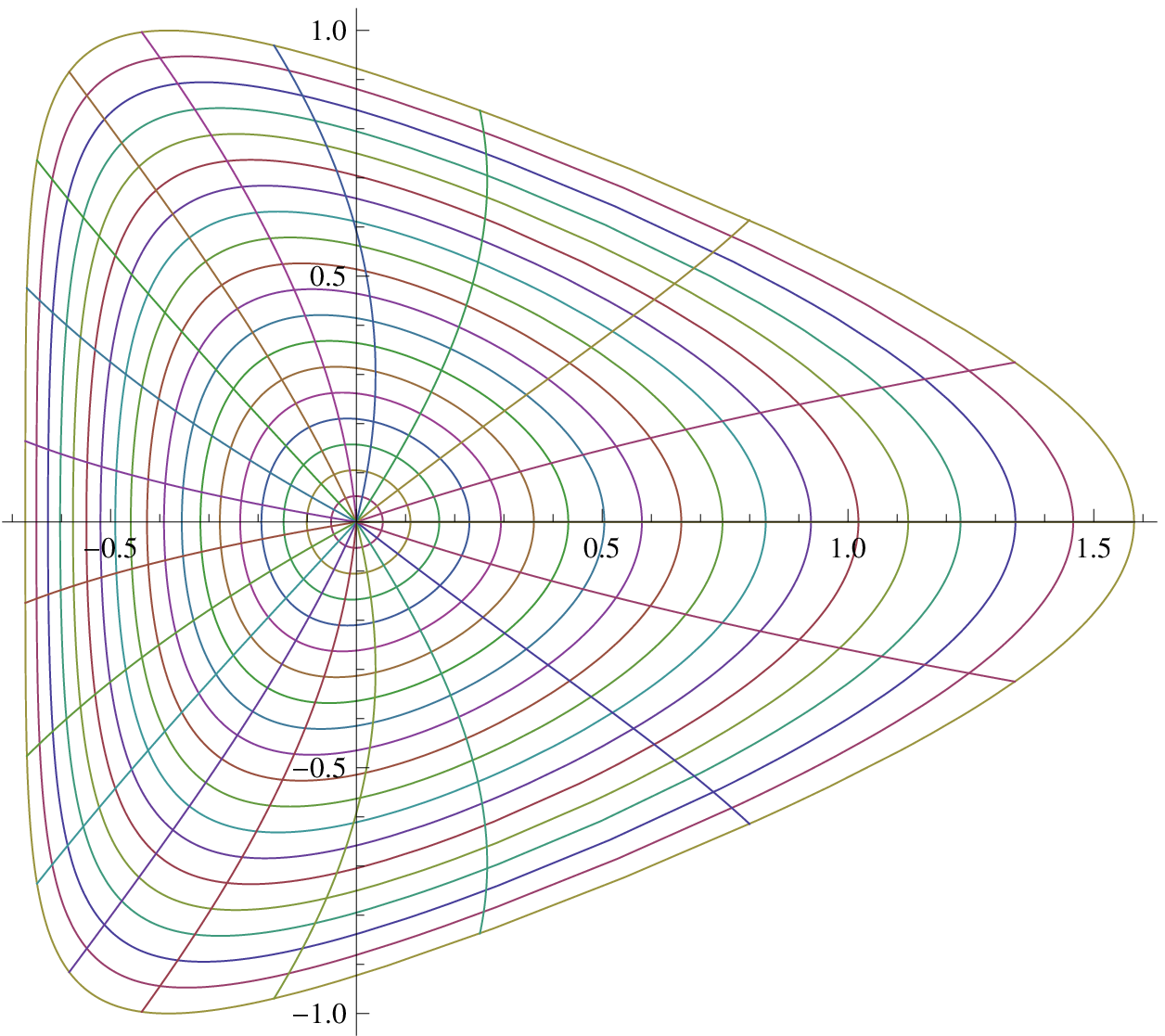}
     \caption{The range  $f_{1+\sqrt{2},0}(\ID)$}
     \label{PS5_Fig1a}
   \end{subfigure}
   ~~~
    \begin{subfigure}[b]{0.5\textwidth}
      \centering
      \includegraphics[height=6cm, width=5.5cm]{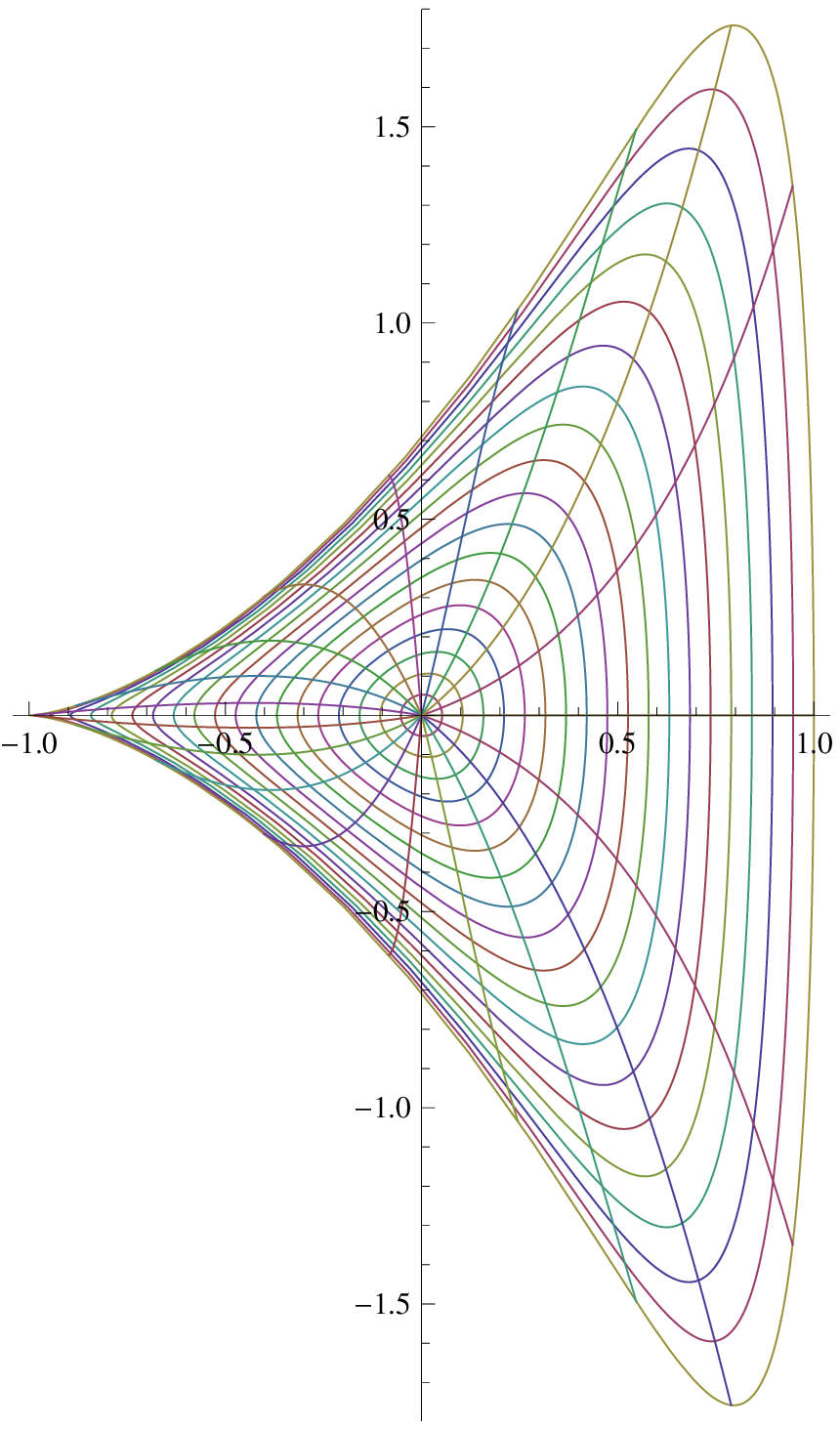}
      \caption{The range $F_{1+\sqrt{2},0}(\ID)$}
      \label{PS5_Fig1b}
   \end{subfigure}\\
   \begin{subfigure}[b]{0.5\textwidth}
     \centering
     \vspace{0.1cm}
     \includegraphics[height=6cm, width=5.5cm]{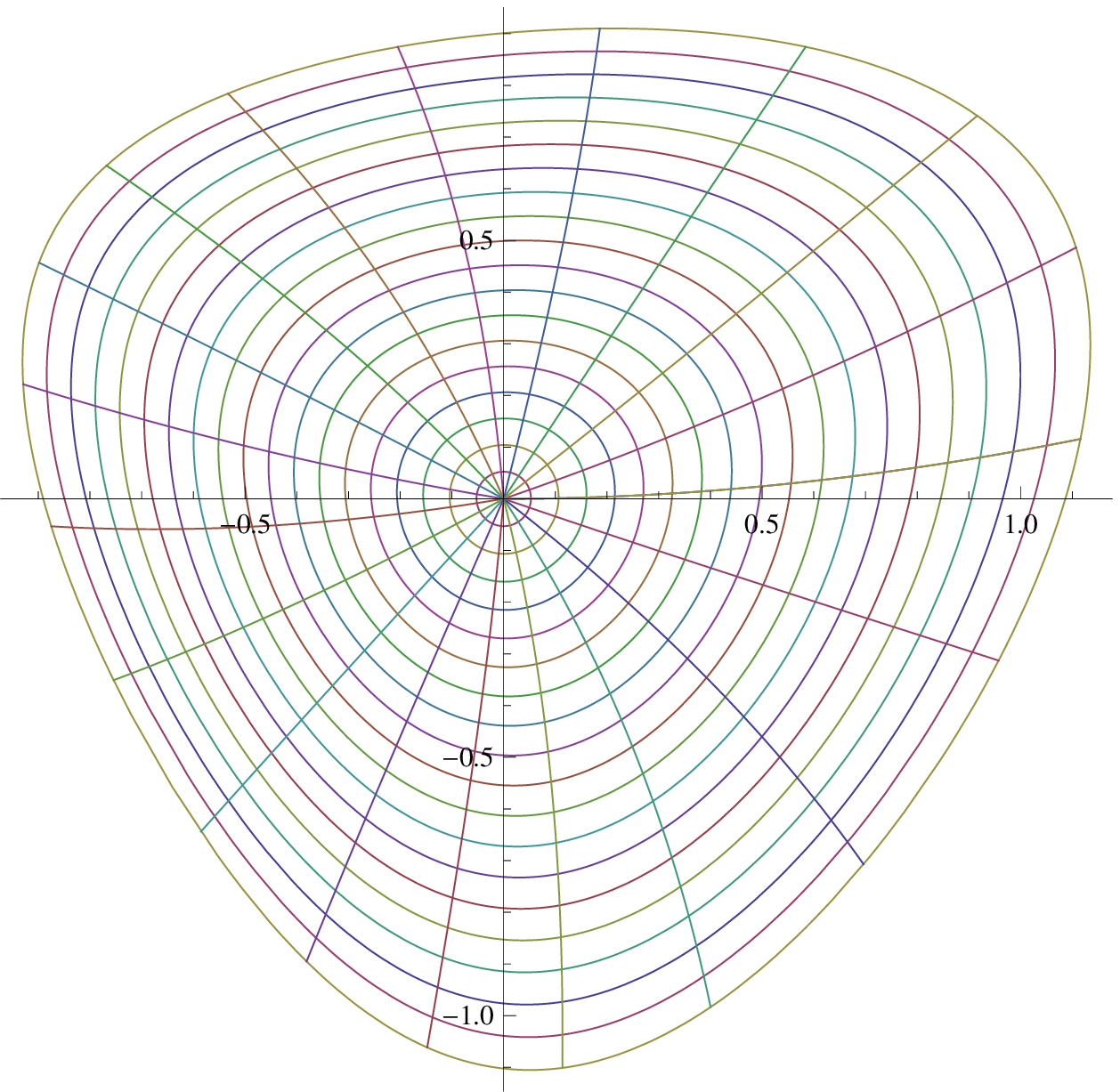}
     \caption{The range $f_{5,i}(\ID)$}
     \label{PS5_Fig1e}
   \end{subfigure}
   ~~~
   \begin{subfigure}[b]{0.5\textwidth}
      \centering
      \includegraphics[height=6cm, width=5.5cm]{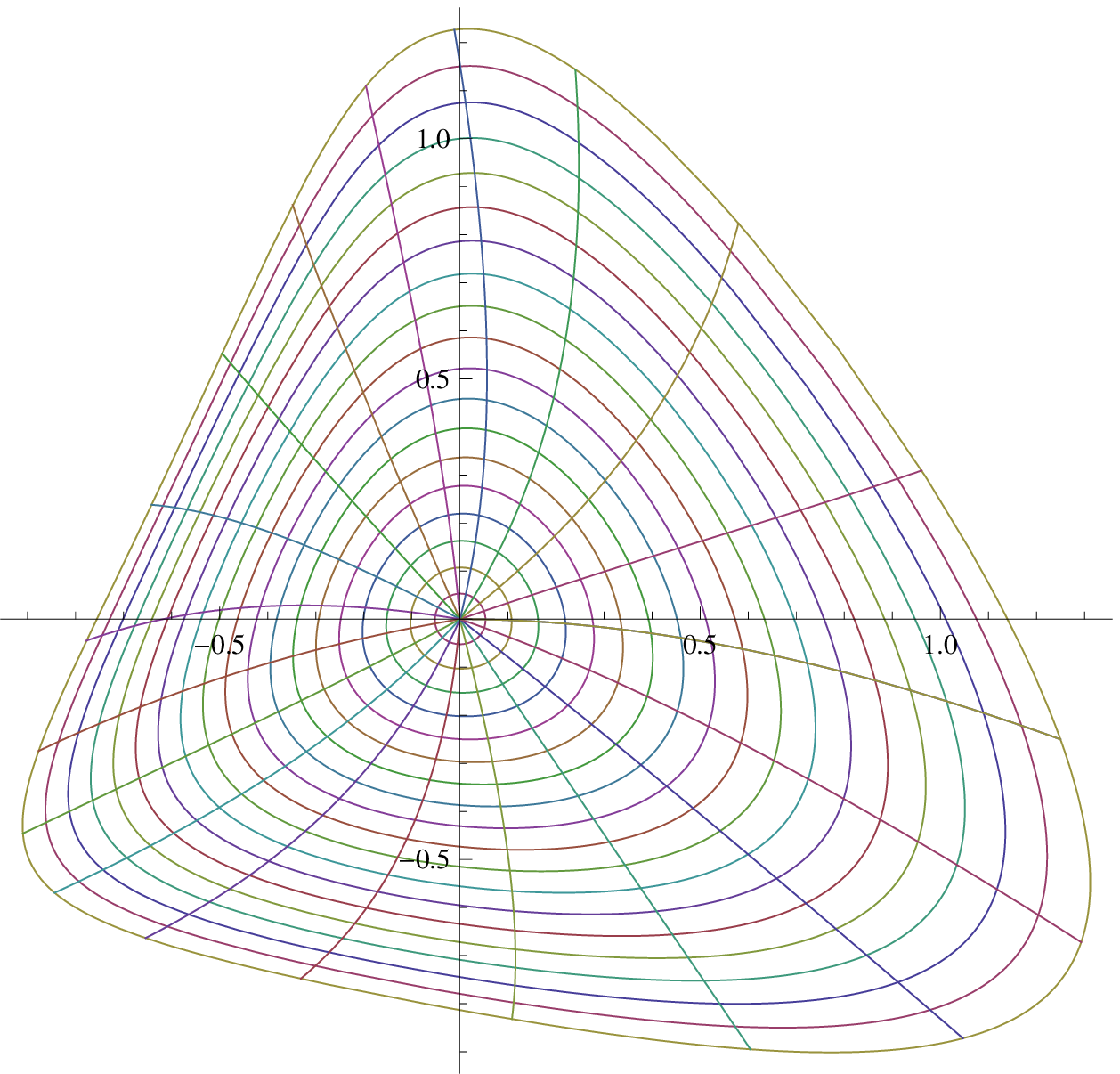}
      \caption{The range  $F_{5,i}(\ID)$}
      \label{PS5_Fig1f}
   \end{subfigure}

\caption{The images of unit disk under $f_{a, \lambda}(z)$ and $F_{a, \lambda}(z)$ for certain values of $a$ and $\lambda$.}
\label{PS5_Fig2}
\end{figure}

Since $f_{a, \lambda}$ defined by \eqref{faz} is {\bf SHC} for all real number $a$ such that $|a| \geq 1+\sqrt{2}$,
we can use Theorem \Ref{SH_Alexander} to construct functions to belong to the class ${\bf SHS}^*$. Define
$$H(z)= \frac{a z}{a-z} , ~\mbox{ and }~ G(z)=\frac{-z^2}{a-z}.
$$
From Theorem \Ref{SH_Alexander}, it is clear that
$$F_{a, \lambda}(z)=\frac{a z}{a-z} -\lambda \overline{\left (\frac{z^2}{a-z}\right )}
$$
is univalent and starlike in $\ID$ for all $\lambda$ such that $|\lambda|=1$ and for any real number $a$ such that $|a| \geq 1+\sqrt{2}$.
We observe that $F_{a, \lambda}(z)=-F_{-a, \lambda}(-z)$. From the definition of
$f_{a, \lambda}(z)$ and $F_{a, \lambda}(z)$, it is a simple exercise to see that
$$\lim_{|a| \rightarrow \infty} f_{a, \lambda}(z) = \lim_{|a| \rightarrow \infty} F_{a, \lambda}(z)=z.
$$
The images of the unit disk $\ID$ under $f_{a, \lambda}(z)$ and $F_{a, \lambda}(z)$ for certain values of
$a$ and $\lambda$ are shown in Figure \ref{PS5_Fig2} as plots of the images of equally spaced radial segments and concentric circles.

\section{Linear and Affine Invariant Families of Harmonic Mappings}\label{PS5sec3}
The class $\mathcal{S}_H(\mathcal{S})$ has several special properties.  For instance, if $f \in \mathcal{S}_H(\mathcal{S})$,
then the function $(f + c \overline{f})/(1+cb_1) \in \mathcal{S}_H(\mathcal{S})$ for all $c \in \mathbb{D}$. This property is
called as affine invariance. Similarly, if
$f \in \mathcal{S}_H(\mathcal{S})$, then for each $\zeta \in \mathbb{D}$, the function $F$ defined by
$$ F(z) = \frac{f(\frac{z+\zeta}{1+\overline{\zeta}z})-f(\zeta)}{(1-|\zeta|^2)h'(\zeta)}
$$
belongs to the class $\mathcal{S}_H(\mathcal{S})$. This is called as linear invariance property \cite{Sheil-Small}. Thus, the family
$\mathcal{S}_H(\mathcal{S})$ is an affine and linear invariant family. Many interesting results have been proved in
the literature for different classes of linear and affine invariant family of harmonic functions. In fact the family
$\mathcal{S}_H$ is invariant under normalized affine and linear transforms. The growth theorem
(see \cite[p.97, Theorem]{Duren:Harmonic}) and the covering theorem for $\mathcal{S}^0_H$ may now be recalled.

\begin{Thm}\label{Theorem_Growth}
Let $\alpha$ be the supremum of $|a_2|$ among all functions $f \in \mathcal{S}_H$. Then, every function
$f \in \mathcal{S}^0_H$ satisfies the inequalities
\be\label{Growth_Ineq}
\frac{1}{2\alpha}\left[1-\left(\frac{1-r}{1+r}\right)^\alpha \right] \leq |f(z)| \leq
\frac{1}{2\alpha}\left[\left(\frac{1+r}{1-r}\right)^\alpha -1 \right], ~~~ r = |z| < 1. \nonumber
\ee
In particular, the range of each function $f \in \mathcal{S}^0_H$ contains the disk $|w| < \frac{1}{2\alpha}$.
\end{Thm}

Let $\mathbb{L}$ be an arbitrary family of locally univalent harmonic functions $f=h+\overline{g}\in {\mathcal H}$, where $h$ and $g$
have the form \eqref{PSSerRep}, such that $\mathbb{L}$ is closed under normalized affine and linear transformations and
$\mathbb{L}^0 = \left\{f \in \mathbb{L}:\,b_1
= 0 \right\}$. Let $\alpha_0$ and $\beta_0$ be the supremum of $|a_2|$ and $|b_2|$, respectively, among all functions
$f \in \mathbb{L}^0$ of the form \eqref{PSSerRep} with $b_1=0$. In \cite{Graf} and \cite{GrafEye}, the authors studied the classes $\mathbb{L}$ and $\mathbb{L}^0$ and derived
the following interesting results.

\begin{Thm}\label{Thm_Jacob}
Let $f \in \mathbb{L}$ with $b_1=f_{\overline{z}}(0)$. Then,
\begin{itemize}
\item[(a)] the Jacobian $J_f$ of the mapping $f$ with any $z \in \mathbb{D}$ satisfies the
bounds
\be\label{Jacob_Ineq}
(1-|b_1|^2)\frac{(1-|z|)^{2\alpha_0 - 2}}{(1+|z|)^{2\alpha_0 + 2}} \leq J_f(z) \leq
(1-|b_1|^2)\frac{(1+|z|)^{2\alpha_0 - 2}}{(1-|z|)^{2\alpha_0 + 2}}. \nonumber
\ee
\item[(b)] for any $z$ with $0 < |z|=r < 1$, the inequalities
$$
|h'(z)| \leq (1+r|b_1|)\frac{(1+r)^{\alpha_0 - 3/2}}{(1-r)^{\alpha_0 + 3/2}}, ~\mbox{ and }~
|g'(z)| \leq (r+|b_1|)\frac{(1+r)^{\alpha_0 - 3/2}}{(1-r)^{\alpha_0 + 3/2}}
$$
hold.
\end{itemize}
These bounds are sharp for the class of univalent close-to-convex harmonic functions. The equality is attained for the close-to-convex functions
$f(z)=K(z)+\overline{b_1K(z)}$, where $K(z)$ is the harmonic Koebe function.
\end{Thm}

\begin{Thm}\label{Thm_Curv}
Let $f \in \mathbb{L}$ with $b_1=f_{\overline{z}}(0)$. Then for $z$ with $0 < |z|=r < 1$, the following bounds
for the curvature $k_f(z)$ of the image of the circle $\{|\zeta|=r\}$ under the mapping $f$ are valid:
\beq\label{Curv_Ineq}
k_f(z) &\leq& \frac{(1+|b_1|)}{(1-|b_1|)^2} \left( \frac{1+r}{1-r} \right)^{\alpha_0 + 3/2}
\frac{r^2+2r(\alpha_0+\beta_0)+1}{r}, \nonumber \\
k_f(z) &\geq& \frac{(1-|b_1|)}{(1+|b_1|)^2} \left( \frac{1-r}{1+r} \right)^{\alpha_0 + 3/2}
\frac{r^2-2r(\alpha_0+\beta_0)+1}{r} ~~\mbox{ if }~~ 0 < r \leq \rho, ~~and \nonumber \\
k_f(z) &\geq& \frac{(1+|b_1|)}{(1-|b_1|)^2} \left( \frac{1+r}{1-r} \right)^{\alpha_0 + 3/2}
\frac{r^2-2r(\alpha_0+\beta_0)+1}{r} ~~\mbox{ if }~~ \rho < r < 1, \nonumber
\eeq
where $\rho = \alpha_0+\beta_0 - \sqrt{(\alpha_0+\beta_0)^2 -1}$.
The inequality $|z| \leq \rho$ determines the maximal disk, where any function $f \in \mathbb{L}$ is convex and univalent.
\end{Thm}

\section{Applications of Theorem \ref{PS5th1}}\label{PS5sec4}
We have already proved the sharp coefficient bounds for the classes $\mathcal{S}^0_H(\mathcal{S})$.
As an application of Theorem \ref{PS5th1} we can now prove sharp coefficient bounds for the
class $\mathcal{S}_H(\mathcal{S})$, and using Theorems \Ref{Theorem_Growth}, \Ref{Thm_Jacob},
and \Ref{Thm_Curv}, we can derive interesting results for the classes $\mathcal{S}^0_H(\mathcal{S})$
and $\mathcal{S}_H(\mathcal{S})$.

As mentioned in the book of Duren (see \cite{Duren:Harmonic}), the only property of the class
$\mathcal{S}_H$ essential to the proof of Theorem \Ref{Theorem_Growth} is its affine and linear invariance. Moreover,
the theorem remains valid for any subclass of $\mathcal{S}_H$ that is invariant under normalized affine and linear
transformation. We have already pointed out that the class $\mathcal{S}_H(\mathcal{S})$ is an affine and linear
invariant family. Hence, Theorem \Ref{Theorem_Growth} is applicable to the class $\mathcal{S}_H(\mathcal{S})$.
Replacing $\mathcal{S}_H$ by $\mathcal{S}_H(\mathcal{S})$ in Theorem \Ref{Theorem_Growth} and applying
Corollary \ref{PS5Cor1} we get the following result. So we omit the details of the proof of these theorems here.

\begin{thm}\label{Growth_sub}
Every function $f \in \mathcal{S}^0_H(\mathcal{S})$ satisfies the inequalities
\be\label{Growth_ineq_sub}
\frac{1}{6}\left[1-\left(\frac{1-r}{1+r}\right)^3 \right] \leq |f(z)| \leq
\frac{1}{6}\left[\left(\frac{1+r}{1-r}\right)^3 -1 \right], ~~~ r = |z| < 1. \nonumber
\ee
In particular, the range of each function $f \in \mathcal{S}^0_H(\mathcal{S})$ contains the disk $|w| < \frac{1}{6}$.
The above inequalities are sharp and the equality is attained for the harmonic Koebe function $K(z)$ and its rotations.
\end{thm}

By taking $\mathbb{L} = \mathcal{S}_H(\mathcal{S})$ in Theorems \Ref{Thm_Jacob}, and \Ref{Thm_Curv}
and applying Theorem \ref{PS5th1} we get the following results.

\begin{thm}\label{thm_jacob_sub}
Let $f \in \mathcal{S}_H(\mathcal{S})$ with $b_1=f_{\overline{z}}(0)$. Then the Jacobian $J_f$ of the mapping $f$
with any $z \in \mathbb{D}$ satisfies the bounds
\be\label{jacob_ineq}
(1-|b_1|^2)\frac{(1-|z|)^3}{(1+|z|)^7} \leq J_f(z) \leq
(1-|b_1|^2)\frac{(1+|z|)^3}{(1-|z|)^7}. \nonumber
\ee
These bounds are sharp. The equality is attained for the close-to-convex functions $f(z)=K(z)+\overline{b_1K(z)}$.
\end{thm}


\begin{thm}\label{thm_dist_sub}
Let $f \in \mathcal{S}_H(\mathcal{S})$ with $b_1=f_{\overline{z}}(0)$. Then for any $z$ with $0 < |z|=r < 1$ the inequalities
$$|h'(z)| \leq (1+r|b_1|)\frac{(1+r)}{(1-r)^4} ~\mbox{ and }~
|g'(z)| \leq (r+|b_1|)\frac{(1+r)}{(1-r)^4}
$$
hold. These bounds are sharp. The equality is attained for the close-to-convex functions
$f(z)=K(z)+\overline{b_1K(z)}$.
\end{thm}

\begin{thm}\label{thm_curv_sub}
Let $f \in \mathcal{S}_H(\mathcal{S})$ with $b_1=f_{\overline{z}}(0)$. Then for any $z$ with $0 < |z|=r < 1$ the following bounds
for the curvature $k_f(z)$ of the image of the circle $\{|\zeta|=r\}$ under the mapping $f$ are valid:
\beq\label{curv_ineq_sub}
k_f(z) &\leq& \frac{(1+|b_1|)}{(1-|b_1|)^2} \left( \frac{1+r}{1-r} \right)^4
\frac{r^2+6r+1}{r}, \nonumber\\
k_f(z) &\geq& \frac{(1-|b_1|)}{(1+|b_1|)^2} \left( \frac{1-r}{1+r} \right)^4
\frac{r^2-6r+1}{r} ~~\mbox{ if }~~ 0 < r \leq \rho, \nonumber ~~and \\
k_f(z) &\geq& \frac{(1+|b_1|)}{(1-|b_1|)^2} \left( \frac{1+r}{1-r} \right)^4
\frac{r^2-6r+1}{r} ~~\mbox{ if }~~ \rho < r < 1, \nonumber
\eeq
where $\rho = 3 - 2 \sqrt{2}$. Moreover, every function $f \in \mathcal{S}_H(\mathcal{S})$ maps the disk
$|z| < 3 - 2 \sqrt{2}\approx0.171572875$ onto a convex domain. That is, for $f \in \mathcal{S}_H(\mathcal{S})$,
the radius of convexity is $3 - 2 \sqrt{2}$.
\end{thm}

We remark that the number $3 - 2 \sqrt{2}$ is the conjectured value by Sheil-Small \cite{Sheil-Small} for
the radius of convexity of $f \in \mathcal{S}_H$. 
%
%


\subsection*{Acknowledgements}
The second author thanks Council of Scientific
and Industrial Research (CSIR), India, for providing financial support in the form of a Senior Research
Fellowship to carry out this research.

\end{document}